\documentclass{amsproc}

\newtheorem{theorem}{Theorem}[section]

\theoremstyle{definition}

\theoremstyle{remark}
\newtheorem{remark}[theorem]{Remark}

\numberwithin{equation}{section}
\def\R{{\mathbb{R}}}
\def\N{{\mathbb{N}}}
\def\C{{\mathbb{C}}}
\def\E{{\mathcal{E}}}

\begin{document}

\title {Asymptotic parabolicity for strongly damped wave equations}

\author{Genni Fragnelli}
\address{Dipartimento di Matematica\\ Universit\`{a} di Bari "Aldo Moro"\\
Via E. Orabona 4\\ 70125 Bari - Italy}
\email{genni.fragnelli@uniba.it}

\author{ Gis\`{e}le Ruiz Goldstein}
\address{Department of Mathematical Sciences\\ University of Memphis\\
Memphis, TN 38152, USA} \email{ggoldste@memphis.edu}

\author{
Jerome A. Goldstein }
\address{
Department of Mathematical Sciences\\ University of Memphis\\
Memphis, TN 38152, USA } \email{jgoldste@memphis.edu}

\author{Silvia Romanelli}
\address{Dipartimento di Matematica\\ Universit\`{a} di Bari "Aldo Moro"\\
Via E. Orabona 4\\ 70125 Bari - Italy} \email{
silvia.romanelli@uniba.it}

\date{September 11, 2012 and, in revised form, ........}

\dedicatory{This paper is dedicated to Fritz Gesztesy on his 60th
birthday .}


\begin{abstract}
For $S$ a positive selfadjoint operator on a Hilbert space,
\[
\frac{d^2u}{dt}(t) + 2 F(S)\frac{du}{dt}(t) + S^2u(t)=0
\]
describes a class of wave equations with strong friction or damping
if $F$ is a positive Borel function. Under suitable hypotheses, it
is shown that
\[
u(t)=v(t)+ w(t)
\]
where $v$ satisfies
\[
2F(S)\frac{dv}{dt}(t)+ S^2v(t)=0
\]
and
\[
\frac{w(t)}{\|v(t)\|} \rightarrow 0, \;  \text{as} \; t \rightarrow
+\infty.
\]
The required initial condition $v(0)$ is given in a canonical way in
terms of $u(0)$, $u'(0)$.
\end{abstract}
\maketitle

\section{Introduction}
Let $S$ be an injective nonnegative selfadjoint operator on a
complex Hilbert space $\mathcal{H}$. That is $S=S^* \ge 0$, $0 \not
\in \sigma_p(S).$ Consider the damped wave equation
\begin{equation}\label{1.1}
u''(t)+2Bu'(t) +S^2u(t)=0, \quad t \ge 0,
\end{equation}
with initial conditions
\begin{equation}\label{1.2}
u(0)=f, \quad u'(0)=g;
\end{equation}
here $'=d/dt$. When $B=0$, \eqref{1.1} reduces to the wave equation,
and the corresponding heat equation normally considered is
\[
v'(t)+S^2v(t)=0.
\]
In this paper we take $B$ to be a positive selfadjoint operator
which commutes with $S$ and is "smaller than $S$". More precisely,
we assume that
\begin{equation}\label{1.3}
0 = \inf \sigma (S),
\end{equation}
i.e., $0$ is the spectrum of $S=S^*\geq 0,$ but is not an
eigenvalue, $F$ is a continuous function from $(0, +\infty)$ to $(0,
+\infty)$, $B=F(S)$, and $F$ satisfies: there exists $\gamma > 0$
such that
\begin{equation}\label{1.4}
\begin{cases}
F(x) > x  &\text{for} \; 0< x<\gamma,\\
F(\gamma)= \gamma,\\
F(x) <x  &\text{for} \; x> \gamma,\\
\limsup_{x \rightarrow 0^+} F(x) < +\infty,\\
\liminf_{x \rightarrow +\infty} ((1- \delta)x-F(x))\ge 0, & \text{
                               for some} \; \delta > 0.

\end{cases}
\end{equation}

We make a further comment on \eqref{1.3}. Think of $S^{2}$ as
$-\Delta $ with suitable boundary conditions, acting on \
$H=L^{2}(\Omega )$ where $\Omega $ is a domain in $\mathbb{R}^{N}.$
Then \eqref{1.3} implies that $S$ cannot have compact resolvent, and
we are thus led to work exclusively in unbounded domains $\Omega .$

The operator $B$ represents a general friction coefficient. The most
common case is the telegraph equation in which case
\[
B= aI
\]
where $a$ is a positive constant. In this case \eqref{1.4} holds
with $\gamma =a$. Another simple case is
\[
B= a S^\alpha
\]
where the constants $a$, $\alpha$ satisfy
\[
a >0, \quad \alpha \in [0,1).
\]
In this case, $\gamma= a^{\frac{1}{1-\alpha}}$ in \eqref{1.4}. The
only interesting case is when $S$ is unbounded. The {\it strongly
damped wave equation} refers to the case when $B$ is also unbounded.

Our main result, Theorem \ref{theorem3.1}, can be stated informally
as follows. Let $S, B=F(S)$ be as above and suppose $f, g$ are such
that \eqref{1.1}, \eqref{1.2} has a unique solution $u$. Consider
the corresponding first order (in $t$) equation, obtained by erasing
$u''(t)$ in \eqref{1.1} and replacing $u$ by $v$:
\begin{equation}\label{1.5}
2Bv'(t)+ S^2v(t)=0, \quad t \ge 0,
\end{equation}
with initial condition
\begin{equation}\label{1.6}
v(0)=h.
\end{equation}
This vector $h$ is given by
\[
h= \frac{1}{2} \chi_{(0,
\gamma)}(S)\{(F(S)^2-S^2)^{\frac{1}{2}}(F(S)f+g)+f\},
\]
a formula which will be derived and explained later; and our
conclusion requires that $f,g$ are such that $ h \neq 0$. We will
show that
\begin{equation}\label{1.7}
l(t):= \frac{\| u(t)- v(t)\|}{\|v(t)\|} \rightarrow 0,
\end{equation}
as $t \rightarrow 0$, and we will find closed subspaces $\E_n$ of
$\mathcal{H}$ such that $\E_n \subset \E_{n+1}$,
\[
\bigcup_{n=1}^\infty\E_n \; \text{is dense in} \; \mathcal {H},
\]
and
\[l(t) \le C_n e^{-\epsilon_nt}
\]
for $f,g \in \E_n$ where $C_n, \epsilon_n$ are positive constants.
In general, there is no rate of convergence in \eqref{1.2} that
works for all solutions.

The point of the theorem  is that, for large times the solution of
the "hyperbolic equation" \eqref{1.1} looks like the solution of the
"parabolic equation" \eqref{1.5}. In the telegraph equation case
when $B= F(S) = aI$, \eqref{1.5} becomes
\[
2av'(t)+ S^2v(t)=0,
\]
which, with \eqref{1.6}, is solved by
\[
v(t)= e^{-\frac{t}{2a}S^2}h.
\]
In the case of strong damping, the solution of the limiting
parabolic problem is
\[
v(t)= e^{-\frac{t}{2} B^{-1}S^2}h.
\]

Think of $S= (-\Delta)^{\frac{1}{2}}$ on $L^2(\R^N)$ and
\[
B= a S^\alpha= a (-\Delta)^{\frac{\alpha}{2}},
\]
$0<\alpha<1.$ Then $B^{-1}S^2 = \displaystyle \frac{1}{a}
(-\Delta)^{1-\frac{\alpha}{2}}$ with domain
\[
D(B^{-1}S^2) = H^{2-\alpha}(\R^N),
\]
while
\[
D(S^2) = H^2(\R^N);
\]
here we use the standard Sobolev space notation. Thus $B^{-1}S^2$ is
a pseudodifferential operator of lower order $2-\alpha$ than that of
the Laplacian, unless $\alpha=0$ in which case we have the telegraph
equation.

It has long been "known" that the telegraph equation is an
asymptotic approximant for the heat equation, especially in the case
of $S^2= - d^2/dx^2$ on $L^2(\R)$. The pioneers in this area include
G. I. Taylor \cite{t} in $1922$, S. Goldstein \cite{gs} in $1938$
and M. Kac \cite{ka} in $1956$.

Some of the associated random walk ideas are discussed in
\cite{egl}, which eventually led to \cite{ceg}, the main theorem in
which is the special case of our main theorem here with $B=aI$. The
importance of the interesting case of strong damping was recognized
by Fritz Gesztesy and is discussed in detail in \cite{gght}. For
additional results on strong damping, see \cite{gr}.

Section \ref{section2} reviews some spectral theory. Section
\ref{section3} contains the statement and the proof of our main
result. The proof is patterned after that in \cite{ceg}, but it in
fact is simplified and streamlined. Section \ref{section4} contains
examples.

\section{Selfadjoint and normal operators}\label{section2}
Let $S$ be a selfadjoint operator on $\mathcal {H}$ with spectrum
$\sigma (S)$. By the spectral Theorem there exists an $L^2$ space
$L^2(\Lambda,\Sigma,\nu)$ and a unitary operator
\[
W: \mathcal {H}\rightarrow L^2(\Lambda,\Sigma,\nu)
\]
such that $S$ is unitarily equivalent, via $W$, to the maximally
defined operator of multiplication by a $\Sigma-$measurable function
\[
m: \Lambda \rightarrow \sigma(S)\subset \R,
\]
i.e.,
\[
Sf= W^{-1} M_m Wf,
\]
for
\[
f \in D(S) =\{W^{-1} g \in \mathcal {H}: \; mg \in L^2(\Lambda,
\Sigma,\nu)\}
\]
and
\[
(M_m g)(x) = m(x)g(x), \quad \text{for} \; x \in \Lambda, g \in
L^2(\Lambda, \Sigma, \nu).
\]

Two selfadjoint operators $S_1$, $S_2$ {\it commute} if and only if
the bounded operators
\[
(\lambda_1 I-S_1)^{-1}, (\lambda_2I- S_2)^{-1}
\]
commute for all $\lambda_1, \lambda_2 \in \C\setminus \R$ if and
only if
\[
e^{itS_1} \quad \text{and} \quad e^{isS_2}
\]
commute for all $t, s \in \R$. Similarly, two normal operators $N_1,
N_2$ with
\[
\sup \text{Re} \;\sigma (N_j) < +\infty, \quad j = 1,2
\]
commute if and only if $e^{tN_1}$ and $e^{sN_2}$ commute for all $t,
s \ge 0$; here $N$ is  normal means $N= S_1+ iS_2$ where $S_1$,
$S_2$ are commuting selfadjoint operators.

The functional calculus for $S$ selfadjoint says that for every
Borel measurable function $F$ from $\sigma (S) \subset \R$ to $\C$,
$F(S)$ defined by
\[
F(S)= W^{-1}M_{F(m)}W
\]
is normal, and any two of these operators commute. Moreover,
\[
F\rightarrow F(S)
\]
is linear and is an algebra homomorphism, thus
\[
F_1(S)F_2(S)= (F_1F_2)(S),
\]
etc. Also, $F(S)$ is bounded on $\mathcal H$ if and only if $F$ is
 bounded on $\sigma (S)$, and $F(S)$ is selfadjoint if and only if
$F$ is real valued. And for $S= S^*$, $F(S)$ is semibounded (above
or below) if and only if $F(\sigma (S))$ is, in $\R$.

In particular, for $\Gamma$ a Borel set in $[0, +\infty)$, $P_\Gamma
= \chi_\Gamma(S)$ is the orthogonal projection of $\mathcal H$ onto
$\chi_\Gamma(S)(\mathcal H)$; and
\[
P_\Gamma F (S) =F(S)P_\Gamma = P_\Gamma F(SP_\Gamma)
\]
is the part of $F(S)$ in $\Gamma$, and its spectrum is contained in
$\Gamma$.

If $F_1, F_2$ are complex Borel functions on $\sigma(S)$ that are
bounded above, it follows that $F_j(S)$ and $\Sigma_{k=1}^nF_k(S)$
generate ($C_0$) semigroups on $\mathcal H$ and
\begin{equation}\label{2.1}
e^{t\Sigma_{k=1}^nF_k(S)} = \Pi_{k=1}^ne^{tF_k(S)},
\end{equation}
and the product can be taken in any order. Finally, if $L= F(S) =L^*
\ge 0$, then $[L]^{\frac{1}{2}}$ denotes the unique nonnegative
square root of $L$.

For more on the spectral theorem, see the books \cite{g}, \cite{k},
\cite{l}.
\section{The main result}\label{section3}
Consider the problem \eqref{1.1}, \eqref{1.2}, which we rewrite as
\begin{equation}\label{3.1}
u'' + 2Bu'+ S^2u =0, \quad t \ge 0,
\end{equation}
\begin{equation}\label{3.2}
u(0)=f, u'(0)=g,
\end{equation}
where $S=S^* \ge 0$ on $\mathcal H$,
\begin{equation}\label{3.3}
\inf \sigma (S)= 0 \not \in \sigma_p(S),
\end{equation}
$B= F(S)$ where $F$ is a continuous function from $(0, +\infty)$ to
$(0, +\infty)$ which is bounded near zero and strictly less than the
identity function near infinity, in the sense that for some $\delta
> 0$,
\begin{equation}\label{3.4}
\liminf_{x\rightarrow +\infty} ((1-\delta) x- F(x)) \ge0.
\end{equation}
We also assume there exists $\gamma > 0$ such that \eqref{1.4}
holds, namely
\begin{equation}\label{3.5}
\begin{cases}
F(x) > x  &\text{for} \; 0< x<\gamma,\\
F(\gamma)= \gamma,\\
F(x) <x  &\text{for} \; x> \gamma,\\
\limsup_{x \rightarrow 0^+} F(x) < +\infty.\\
\end{cases}
\end{equation}
Let $\Gamma$ be the open interval $(0, \gamma)$ and let
\begin{equation}\label{3.6}
P_\Gamma = \chi_\Gamma (S).
\end{equation}

\medskip
\begin{theorem}\label{theorem3.1}Assume all the statements in the
above paragraph. Let $v$ be the solution of the corresponding
"parabolic" equation
\begin{equation}\label{3.7}
2Bv'+ S^2v=0,
\end{equation}
obtained by deleting the second derivative term in \eqref{3.1}. Let
\begin{equation}\label{3.8}
v(0)= h:= \frac{1}{2}P_{\Gamma}
(f+[(B^2-S^2)P_{\Gamma}]^{\frac{1}{2}}(Bf+g)).
\end{equation}
Then, for $u$ the solution of \eqref{3.1}, \eqref{3.2},
\begin{equation}\label{3.9}
u(t)= v(t)(1+o(t))
\end{equation}
holds as $t \rightarrow +\infty$, provided $h \neq0$. Moreover, if
$\Gamma_n =\left[\frac{1}{n}, \delta- \frac{1}{n}\right]$ and if $0
\neq h \in P_{\Gamma_n}(\mathcal{H})$ for some $n \in \N$, then
\begin{equation}\label{3.10}
u(t) =v(t)(1+o(e^{-\epsilon_nt}))
\end{equation}
for some $\epsilon_n >0$.
\end{theorem}

\begin{remark} Note that
\[
P_{\Gamma_n}(\mathcal{H}) \subset P_{\Gamma_{n+1}}(\mathcal{H}),
\]
and $\bigcup_{n=1}^\infty P_{\Gamma_n}(\mathcal{H})$ is dense in
$P_\Gamma(\mathcal{H})$.
\end{remark}

\begin{proof}
First recall that the square root of $(B^2-S^2)P_\Gamma$ refers to
the nonnegative square root. We first show that the problem
\eqref{3.1}, \eqref{3.2} is wellposed, by showing that it is
governed by a ($C_0$) contraction semigroup.

We first treat the case of $B=0$ in \eqref{3.1}.

Rewrite \eqref{3.1}, \eqref{3.2} as, for $U=\begin{pmatrix}
Su\\u'\end{pmatrix}$,
\begin{equation}\label{3.11}
\begin{aligned}
U'= \begin{pmatrix} Su' \\u''\end{pmatrix} &= \begin{pmatrix} 0 &S\\
-S &0\end{pmatrix} \begin{pmatrix}Su\\ u'\end{pmatrix} = GU,\\
U(0) &= L= \begin{pmatrix} Sf \\g\end{pmatrix}.
\end{aligned}
\end{equation}
Let $\mathcal K$ be the completion of $D(S)\oplus D(S)$ in the norm
\[
\left\|\begin{pmatrix}m\\n\end{pmatrix}\right\|_{\mathcal K} =
(\|m\|^2 +\|n\|^2)^{\frac{1}{2}},
\]
where $\|\cdot\|$ is the norm on $\mathcal{H}$. Then $G$ defined by
\eqref{3.11} is skewadjoint on $\mathcal K$ and thus generates a
($C_0$) unitary group by Stone's Theorem (cf. \cite{g}, \cite{k},
\cite{l}). We consider this as a semigroup since we are only
concerned with times $t \ge 0$. Next, \eqref{3.1} can be written in
$\mathcal K$ as
\[
U' = (G+P)U
\]
where
\[
P=\begin{pmatrix} 0&0\\ 0& -2B\end{pmatrix}.
\]
Since $B= B^* \ge 0$,
\[
\begin{aligned}
\left\|P\begin{pmatrix} m \\n\end{pmatrix}\right\|_{\mathcal{K}} &=
\|2Bn\| \le (1-\epsilon)\|Sn\| + C_\epsilon\|n\|\\
&\le
(1-\epsilon)\left\|P\begin{pmatrix}m\\n\end{pmatrix}\right\|_{\mathcal{K}}+
C_\epsilon\left\|\begin{pmatrix}
m\\n\end{pmatrix}\right\|_{\mathcal{K}}
\end{aligned}
\]
for some $\epsilon >0$ and a corresponding $C_\epsilon >0$, for all
$n \in D(S) \subset D(B)$, thanks to \eqref{3.4} and the last line
in \eqref{3.5} (or \eqref{1.4}). Namely, write $B= B_1+B_2:= BP_{(0,
M)} +BP_{[M, \infty)}$, where $M$ is such that
\[x \ge F(x) + \delta x,
\]
i.e.
\[
F(x) \le (1-\delta)x
\]
for $x \ge M$ and $F(x)$ is bounded in $[0, M]$. Thus $B_1$ is
bounded, $B_2=B_2^* \ge 0$ and
\[
\|B_2n\| \le (1-\delta)\|Sn\|+\|B_2\|\|n\|
\]
for all $n \in D(S)$. Thus for $N=
\begin{pmatrix}m\\n\end{pmatrix},$
\[
\|PN\|_{\mathcal K}\le (1-\delta) \|GN\|_{\mathcal K} +
M\|N\|_{\mathcal K}
\]
where $\delta >0$ and $M =\|B_2\|.$ It follows (cf. e.g. \cite{g},
 \cite{k}, \cite{l}) that $G+P$ is $m-$dissipative and generates a
($C_0$) contraction semigroup on $\mathcal K$, since $P$ is
obviously dissipative. Then \eqref{3.1}, \eqref{3.2} has a unique
strongly $C^2$ solution (resp. mild solution) if $f \in D(S^2), g
\in D(S)$ (resp. $ f\in D(S), g \in \mathcal{H})$: cf. \cite[Chapter
2, Theorem 7.8]{g}.

We shall express the unique solution using d'Alembert's formula. We
seek solution of the form
\[
u(t)= e^{tC}h
\]
where $C$ is a Borel function of $S$. By \eqref{3.1}, $C$ must
satisfy
\[
C^2 +2BC+S^2 =0.
\]
Formally,
\[
C=C_{\pm}= -B \pm (B^2 -S^2)^{\frac{1}{2}}.
\]
Selfadjoint operators have many square roots, but nonnegative
selfadjoint operators have unique nonnegative square roots. Thus we
uniquely define $C_\pm$ by
\begin{equation}\label{3.12}
C_{\pm}= -B \pm (Q_0 +iQ)
\end{equation}
where
\begin{equation}\label{3.13}
Q_0 =[(B^2-S^2)\chi_{(0,\gamma)}(S)]^{\frac{1}{2}},
Q=[(S^2-B^2)|\chi_{[\gamma,+\infty)}(S)]^{\frac{1}{2}}.
\end{equation}
Thus the solution $u$ of \eqref{3.1}, \eqref{3.2} can be written as
\[
u(t) = e^{tC_+}h_+ + e^{tC_-}h_-,
\]
where $C_\pm$ are defined by \eqref{3.12}, \eqref{3.13}. There are
strong $C^2$ solutions (resp. mild solutions) if and only if $h_\pm
\in D(S^2)$ (resp. $h_\pm \in D(S)$).

Given $f =u(0), g=u'(0)$, we obtain $h_\pm$ by inverting the $2
\times 2$ system
\[
\begin{aligned}
&f = h_+ +h_-\\
&g= C_+h_+ + C_-h_-.
\end{aligned}
\]
An elementary calculation gives
\begin{equation}\label{3.14}
h_- = \frac{1}{2} (f-(Q_0 +iQ)^{-1}(Bf+g))
\end{equation}
\begin{equation}\label{3.15}
h_+= \frac{1}{2}(f+ (Q_0 +iQ)^{-1}(Bf+g)).
\end{equation}
Write
\[
u =u_1+u_2+u_3
\]
where
\[
u_1(t)= e^{tC_+}P_{(0, \gamma)}h_+,
\]
\[u_2(t)=e^{tC_+}P_{[\gamma, +\infty)}h_+,\]
\[
u_3(t) = e^{tC_-}h_-.\] First,
\[
\|u_3(t)\|= \|e^{-itQ}e^{-tQ_0}e^{-tB}h_-\| \le \|e^{-tB}h_-\|
\]
since $e^{-itQ}$ is unitary and $\|e^{-tQ_0}\| \le 1$. Next,
\[
\|u_2(t)\| =\|e^{itQ}e^{-tB}P_{[\gamma, +\infty)} h_+\| \le
\|e^{-tB}h_+\|. \] The next estimate is the key one.

For
\begin{equation}\label{3.16}
h:= P_{(0, \gamma)}(h_+),\end{equation}
\[\|u_1(t)\| =\|e^{t(-B+Q_0)}h_+\|.\]

We know that $h \in P_{(0, \gamma)}(\mathcal H)$: assume
\begin{equation}\label{3.17}
0 \neq h \in P_{[\delta, \gamma -\delta]}(\mathcal H)=: \mathcal
{H}_\delta
\end{equation}
for some $\delta> 0$. Let
\[
Q_{0\delta} =Q_0P_{[\delta, \gamma -\delta]}.
\]
Then, since $F(x) > x$ on $[\delta,\gamma - \delta]$, $F(x) -x \ge
\epsilon$ on $[\delta, \gamma -\delta]$ for some $\epsilon >0$. Thus
\[
Q_{0\delta} \ge \epsilon I.
\]
Consequently
\[
\|u_1(t) \| \ge e^{\epsilon t} \|e^{-tB}h\|.
\]
It follows  that for some constant $C_0$,
\[
\|u_2(t)\|+ \|u_3(t)\| \le C_0 e^{-\epsilon t}\|u_1(t)\|.
\]
Thus
\begin{equation}\label{3.18}
u(t) = u_1(t) (1+O(e^{-\epsilon t})).
\end{equation}
We must show that this holds with $u_1$ replaced by $v$.

The unique solution of \eqref{3.7} is
\[
v(t)= e^{-\frac{t}{2}B^{-1}S^2}h.
\]
Note that $h$ as defined by \eqref{3.8} is $P_\Gamma h_+$  where
$h_+$ is as in \eqref{3.15}. To compare $u_1$ with $v$, we need
Taylor's formula with integral remainder, which for $g \in C^3
[0,l]$ and for some $l
>0$ says that
\[
g(x) = g(0)+ g'(0)x + \frac{g''(0)}{2} x^2 + \frac{1}{2}\int_0^x
(x-y)^2 g'''(y)dy.
\]
Applying this to
\begin{equation}\label{3.19}
g(x) = 1- (1-x)^{\frac{1}{2}}, \quad 0 <x<1,
\end{equation}
 yields
 \begin{equation}\label{3.20}
g(L) f= \frac{1}{2}Lf +\frac{1}{8} L^2 f + Rf
 \end{equation}
 where $R$ is a bounded operator commuting with $L$ and satisfying
 \[
 R=R^*\ge 0.
 \]
 Consequently
 \begin{eqnarray}
\|u_1(t) -v(t) \| &=\|e^{-tB(-I +B^{-1}Q_0)}h -e^{-\frac{t}{2}
B^{-1}S^2}h\|\nonumber\\  &\;\;=
\|e^{-tB[I-(B^{-2}Q_0^2)^{\frac{1}{2}}]}h-e^{-\frac{t}{2}
B^{-1}S^2}h\|\nonumber\\
&\!\!\!\!\!\!\!\!=\|e^{-\frac{t}{2}B^{-1}S^2}
\{e^{-\frac{t}{8}L^2}e^{-tR}-I \}h\| \label{3.21}
 \end{eqnarray}
by \eqref{3.19}, \eqref{3.20} with $L= B^{-2}(B^2 -S^2)P_{(0,
\gamma)} = (I- B^{-2}S^{-2})P_{(0, \gamma)}.$

We have
\[
\zeta_1 I \le R\le \zeta_2 I
\]
on $\mathcal{H}_\delta$ for some constants $0 < \zeta_1< \zeta_2<
+\infty$. Furthermore, we also have
\[
\zeta_3 I \le L\le \zeta_4 I
\]
on $\mathcal{H}_\delta$ for some positive constants $ \zeta_3,
\zeta_4$. It now follows from \eqref{3.20} that
\[
\|u_1(t) -v(t)\| = \|e^{-\frac{t}{2}
B^{-1}S^2}(I-e^{-\frac{t}{8}L}e^{-tR}h)\|
\]
and
\[
\|e^{-\frac{t}{8}L}e^{-tR}h\| \le e^{-t\zeta_5}\|h\|
\]
where
\[
\zeta_5 =\frac{\zeta_3}{8} + \zeta_1 >0.
\]
Consequently
\[
\|u(t) -v(t)\|\le \|v(t)\| O(e^{-t\zeta_5}).
\]
Combining this inequality with \eqref{3.18} yields the desired
asymptotic relation
\[
\frac{\|u(t)-v(t)\|}{\|v(t)\|} = o(e^{-t\epsilon_\delta})
\]
for some $\epsilon_\delta >0$.

Now let $0 \neq h \in P_{(0, \gamma)} (\mathcal H).$ We must show
that
\begin{equation}\label{3.22}
\frac{\|u(t)-v(t)\|}{\|v(t)\|} \rightarrow 0,\quad \text{and }  t
\rightarrow +\infty.
\end{equation}

We proceed by contradiction. Suppose \eqref{3.22} fails to hold for
some $h \neq 0$ in $P_{(0, \gamma)}(\mathcal H)$. Then, there exists
$\epsilon_1 >0$ and $t_n \rightarrow  +\infty$ such that
\begin{equation}\label{3.23}
\frac{\|u(t_n)-v(t_n)\|}{\|v(t_n)\|}\ge \epsilon_1
\end{equation}
for all $n \in \N$. Choose $\delta >0$ and $\tilde h \in
\mathcal{H}_\delta = P_{[\delta, \gamma-\delta]}(\mathcal H)$
(depending on $\epsilon_1$) such that
\[
\|h -\tilde h\| <\frac{\epsilon_1}{4}
\]
and let $\tilde f$, $\tilde g$ be the corresponding initial data.
Note that
\[
P_{[\gamma, +\infty)}l= P_{[\gamma, +\infty)}\tilde l
\]
for $l =f, g$, and $f$ and $g$ are modified only on the subspace
$P_{\Lambda}(\mathcal H)$
\[
\Lambda :=[\delta -\delta_1, \delta +\delta_1] \cup
[\gamma-\delta-\delta_1, \gamma-\delta+\delta_1],
\]
for some $\delta_1 >0$ which can be chosen to be arbitrarily small.
In particular, given $\epsilon_2>0$ we may choose $\tilde f, \tilde
g$ as above and additionally satisfying
\[
\|f-\tilde f\|+\|g -\tilde g\| <\epsilon_2,
\]
\[
\frac{\|\tilde h\|}{\|h\|} \in [1-\epsilon_2, 1+ \epsilon_2].
\]
It follows that
\[\|u(t) -\tilde u(t) \|, \;\|v(t) -\tilde v(t)\|
<\frac{\epsilon_1}{4}
\]
for all $t >0$. Consequently
\begin{equation}\label{3.24}
\begin{aligned}
\frac{\|u(t) -v(t)\|}{\|v(t)\|} &\le \frac{\|\tilde u(t) -\tilde
v(t)\|}{\|\tilde v(t) \| } \left(
\frac{1+\epsilon_2}{1-\epsilon_2}\right) +\frac{\epsilon_1}{4}\\&
\le \tau_0 e^{-\epsilon_3t} \left(
\frac{1+\epsilon_2}{1-\epsilon_2}\right)
+\frac{\epsilon_1}{4}\rightarrow\frac{\epsilon_1}{4},
\end{aligned}
\end{equation} as $t \rightarrow  +\infty$, since $0\neq \tilde h
\in \mathcal{H}_\delta$, and $\tau_0, \epsilon_3$ are positive
constants depending on $\delta$. But \eqref{3.24} contradicts
\eqref{3.23} for $t=t_n$ with $n$ large enough. It follows that
\eqref{3.22} holds. This completes the proof of Theorem
\ref{theorem3.1}.
\end{proof}
We remark that, in general, there does not exist a rate of
convergence in \eqref{3.22} which works for all $0 \neq h \in
\mathcal{H}_\gamma$. This follows from a careful examination of
\cite{ceg}  and \cite[Section 3]{cggr}.
\section{Examples}\label{section4}
{\bf Example 4.1.} This is the simplest example. Take
\[
S^2 :=- \Delta +w^2 I \quad \text{on}\; L^2(\R^N)
\]
for $w >0$. Similarly we can define
\[
S^2_k := S^{2k}= (- \Delta + w^2I)^k \quad \text{on} \; L^2(\R^N),
\; k \in \N.
\]
Theorem \ref{theorem3.1} applies (with $k \in \N$, $w >0$ fixed) if
we take
\[
B= aS_k^\alpha = a S^{\alpha k}
\]
for $a >0$, $\alpha \in [0,1)$. The resulting damped wave equation
is
\begin{equation}\label{4.1} u_{tt} + 2a(-\Delta
+w^2)^{\frac{\alpha k}{2}} u_t + (-\Delta + w^2)^k u = 0.
\end{equation}
This is a pseudo differential equation unless $\frac{\alpha k}{2}
\in \N_0=\{0,1,2,...\},$ in which case it is a partial differential
equation. The simplest case of this is $k=3$, $\alpha =
\frac{2}{3}$, in which case \eqref{4.1} reduces to
\[
u_{tt} + 2a(-\Delta u_t + w^2 u_t) + (-\Delta+w^2)^3 u=0.
\]
As noted earlier, $\gamma$ is given by
\[
\gamma = a^{\frac{1}{1-\alpha}}
\]
in the case of \eqref{4.1}, for all $k \in \N$.

\medskip
 {\bf Example 4.2.}
We start with a brief summary of the "Wentzell Laplacian" discussed
in \cite{cggr} (cf. also \cite{cfgggor}). Let $\Omega$ be an
unbounded domain in $\R^N$ with nonempty boundary $\partial \Omega$,
such that for every $R>0$ there exists a ball $B(x_R, R)$ in
$\Omega$. Let $A(x)$ be an $N\times N$ matrix for $x \in
\overline{\Omega}$ such that $A(x)$ is real, Hermitian and
\begin{equation}\label{4.2}
\alpha_1|\xi|^2\le A(x) \xi \cdot\xi \le \alpha_2 |\xi|^2
\end{equation}
for all $x \in \overline{\Omega}, \xi \in \R^N$, where
\[
0<\alpha_1\le \alpha_2<\infty
\]
are constants. Similarly let $B(x)$ for $x \in \partial \Omega$ be a
real Hermitian $(N-1) \times (N-1)$ matrix satisfying \eqref{4.2}
with the same $\alpha_1$, $\alpha_2$. Assume $\partial \Omega, A,
B$, and $\gamma, \beta$ (introduced below) are sufficiently smooth.
Define distributional differential operators on $\Omega$ (resp.
$\partial \Omega$) by
\[
Lu_1 =\nabla\cdot (A(x) \nabla u_1),
\]
\[
L_\partial u_2= \nabla_\tau \cdot (B(x)\nabla_\tau u_2)
\]
for $u_1$ (resp. $u_2$) defined on $\Omega$ (resp. $\partial
\Omega$). Here $\nabla_\tau$ is the tangential gradient on $\partial
\Omega$. The wave equation (without damping) we consider is
\begin{equation}\label{4.3}
u_{tt} =Lu \quad \text{in } \Omega,
\end{equation}
\begin{equation}\label{4.4}
Lu + \beta \partial_\nu^A u + \gamma u -q\beta L_\partial u=0 \quad
\text{on } \partial \Omega.
\end{equation}
Here the conormal derivative term is
\[
\partial_\nu^A u =(A\nabla u)\cdot \nu
\]
at $x \in \partial \Omega$, where $\nu$ is the unit outer normal to
$\partial \Omega$ at $x$; $\beta,\gamma\in
C^1(\partial\Omega,\mathbb{R}), \beta>0,\gamma\ge 0$,
 $\beta, \frac{1}{\beta}$ and $\gamma$ are
bounded, and $q \in [0, +\infty)$.

In \cite{cggr}, it is explained how the problem \eqref{4.3},
\eqref{4.4} can be rewritten as
\[
u''+ S^2 u=0
\]
\[u(0)=f, \quad u'(0)=g\]
where the Hilbert space is
\[
\mathcal H = L^2(\Omega, dx)\oplus L^2\left(\partial \Omega,
\frac{d\Gamma}{\beta}\right).
\]
Vectors in $\mathcal H$ are represented by $ u = \begin{pmatrix} u_1
\\u_2\end{pmatrix} $ with $u_1 \in L^2(\Omega, dx)$ and $u_2 \in
L^2\left(\partial \Omega, \frac{d\Gamma}{\beta}\right)$. The norm in
$\mathcal{H}$ is given by
\[
\|u\|_{\mathcal H}=\{\|u_1\|^2_{L^2(\Omega)}
+\|u_2\|^2_{L^2\left(\partial \Omega,
\frac{d\Gamma}{\beta}\right)}\}^{\frac{1}{2}}.
\]
Here we write $d\Gamma $ rather than the usual $dS$ for the element
of surface measure, since the letter $S$ already is being used to
denote the basic operator.

The operator $S^2$ has the matrix representation
\[
S^2 = \begin{pmatrix} -L & 0\\ \beta \partial_\nu^A & \gamma-q\beta
L_\partial\end{pmatrix}.
\]
In \cite{cggr} it was shown that $S=[S^2]^{\frac{1}{2}}$, with a
suitable domain, satisfies
\[
S=S^* \ge 0, \; 0 = \inf \sigma (S), \;0 \not \in \sigma_p(S).
\]
Furthermore, for all $u \in D(S)$, we have $u= \begin{pmatrix}
u_1\\u_2\end{pmatrix}$, where $u_2 = \text{tr} (u_1)$, the trace of
$u_1$. As in Example \ref{4.1}, Theorem \ref{theorem3.1} applies to
\begin{equation}\label{4.5}
u'' + 2aS^{\frac{\alpha k}{2}}u' + S^{2k} u=0,
\end{equation}
\[u(0)=f, \; u'(0)=g, \quad k \in \N.
\]
Again, this is a partial differential equation only when
$\frac{\alpha k}{2} \in \N$. If $k=3$ and $\alpha=\frac{2}{3}$, the
corresponding parabolic problem is
\[
v'+ \frac{1}{2a}S^4 v =0, \quad v(0)=h.
\]
The boundary conditions associated with \eqref{4.5} are
\[ Lw +
\beta\partial_\nu^A w + \gamma w-qL_\partial w = 0 \quad \text {on }
\partial \Omega
\]
for $w =S^{2j} u$, $j=0,1,...,k-1$.

\medskip
 {\bf Example 4.3} The simplest example of unidirectional
waves in one dimension are described by the equation (for $t, x\in
\R$)
\begin{equation}\label{4.6}
u_t = cu_x + bu_{xxx} =: Mu,
\end{equation}
where $(b,c) \in \R^2 \setminus \{(0,0)\}.$ The most common case is
$c \neq  0,b=0$, in which case the corresponding equation for
bidirectional waves is
\[
\left( \frac{\partial}{\partial t} -
M\right)\left(\frac{\partial}{\partial t} + M\right) u = u_{tt}- c^2
u_{xx} = 0.
\]
The case of $b \neq 0$ is the Airy equation, and \eqref{4.6} is the
linearization of the KdV equation
\[
u_t = cu_x + bu_{xxx}+ c_1 u u_x.
\]
For $c= 0 \neq b$, the bidirectional version of \eqref{4.6} is
\[
\left( \frac{\partial}{\partial t} -
M\right)\left(\frac{\partial}{\partial t} + M\right) u = u_{tt}- b^2
u_{xxxxxx} = 0.
\]

Now, let $\mathcal H = L^2(\R)$, $D = \frac{d}{dx}$ and $T= -D^2
=T^*\ge 0$. Let
\[
S^2 = T^3 + a_0 T^2 + a_1T,
\]
where $a_0$, $a_1 \in [0,+\infty)$. Consider
\[
u_{tt}-2au_{xxt}-u_{xxxxxx} + a_0u_{xxxx}-a_1u_{xx}=0,
\]
\[
u(x,0)= f(x), \; u_t(x,0) = g(x).
\]
In this case,
\[
B= aT = F(S) = F(T^3+ a_0T^2+ a_1T).
\]
For $x > 0$, we want to consider the function
\[
G(x) = \frac{1}{a}(x^3+ a_0 x^2 + a_1x),
\]
so that $G(x) -x$ is negative in $(0, \gamma)$ and positive on
$(\gamma,\infty)$ for some $\gamma >0$.  But
\[
\frac{d}{dx}\left(\frac{G(x) -x}{x}\right) = \frac{d}{dx}
\left(\frac{1}{a} (x^2+ a_0x + (a_1-a))\right) = \frac{1}{a} (2x +
a_0) >0,
\]
and $G(x) =x$ for $x \neq 0$ if and only if
\[
x= \frac{1}{2}(-a_0\pm \sqrt{a_0^2-4(a_1-a)}\;).
\]
Thus we get exactly one positive root if and only if
\begin{equation}\label{4.7}
a>a_1+\frac{a_0^2}{4},
\end{equation}
which we assume. It is now elementary to check that $B=F(S)$ and $F$
satisfies the assumptions of Theorem \ref{theorem3.1}. In this case
\[
\gamma =\frac{1}{2}(-a_0 + \sqrt{a_0^2- 4(a_1-a)}\;).
\]

\end{document}